\documentclass[%
 reprint,
 amsmath,amssymb,
 aps,
]{revtex4-2}

\usepackage{graphicx}
\usepackage{dcolumn}
\usepackage{bm}
\usepackage{subcaption}

\usepackage{booktabs}
\usepackage{multirow}

\begin{document}

\preprint{APS/123-QED}

\title{Optimizing Hard Thresholding for Sparse Model Discovery}

\author{Derek W. Jollie}
 \altaffiliation[Also at ]{Department of Physics, Montana State University}
 \email{derek.jollie@student.montana.edu}
\author{Scott G. McCalla}%
 \email{scott.mccalla@montana.edu}
\affiliation{%
 Department of Mathematical Sciences, Montana State University
}%

\date{\today}

\begin{abstract}
Many model selection algorithms rely on sparse dictionary learning to provide interpretable and physics-based governing equations. The optimization algorithms typically use a hard thresholding process to enforce sparse activations in the model coefficients by removing library elements from consideration.  By introducing an annealing scheme that reactivates a fraction of the removed terms with a cooling schedule, we are able to improve the performance of these sparse learning algorithms. 
 We concentrate on two approaches to the optimization, SINDy, and an alternative using hard thresholding pursuit. We see in both cases that annealing can improve model accuracy. The effectiveness of annealing is demonstrated through comparisons on several nonlinear systems pulled from convective flows, excitable systems, and population dynamics. Finally we apply these algorithms to experimental data for projectile motion. 
\end{abstract}

\maketitle

\noindent\emph{Introduction} 
Quantitative science progresses through model discovery, or by constructing mathematical models that reproduce known data. These models created from observed natural phenomena allow researchers to analyze the underlying mechanisms for their importance, forecast future behavior, and seek commonalities among different phenomena from various fields.  Historically model discovery required deep domain expertise, and was especially successful in fields where the underlying physical laws were known, such as Newton's Laws or Quantum Mechanics.  In recent years, Scientific Machine Learning (SciML) has lead to numerous new algorithms combining techniques from sparse 
\cite{bongard2007automated,wang2021bridging,berry2023learning,north2023review,tran2017exact,dam2017sparse,boninsegna2018sparse,messenger2021weak} and deep 
\cite{long2018pde,cai2021physics,cuomo2022scientific,boulle2023mathematical,karniadakis2021physics,lin2023learning,li2020fourier,lu2021learning,liu2023random,jollie2024time,liu2024prose} learning to discover these governing laws directly from data.

More specifically, given time series data $X(t)$ the goal is to find a model $\dot{X}=f(X(t))$ that reproduces $X(t)$.  
Sparse model selection algorithms, such as Sparse Identification of Nonlinear Dynamics (SINDy) \cite{doi:10.1073/pnas.1517384113} and related approaches
\cite{fasel2022ensemble,schaeffer2017sparse,MESSENGER2021110525,schaeffer2017learning,schaeffer2018extracting,fasel2021sindy}, perform exceptionally well on test datasets originating from ordinary and partial differential equations (ODEs and PDEs) in low noise settings.
In these situations, convergence to a local minimizer can be guaranteed \cite{doi:10.1137/18M1189828}. Pairing these algorithms with domain expertise can lead to optimal solutions \cite{viknesh2024adam,MESSENGER2021110525}.
The addition of statistical methods to SINDy also provides a promising approach to strengthening these algorithms such as through ensemble learning \cite{fasel2022ensemble,bekar2024multiphysics}, or \cite{champneys2024bindybayesianidentification,mars2024bayesian,klishin2024statistical} which use a Bayesian framework for SINDy.
 Real data is often limited and exhibits both large measurement and system noise as well as fully stochastic dynamics,
 requiring extra assumptions or techniques to fully discover an interpretable and generalizable model \cite{de2020discovery}.  The current work aims to circumvent some of the issues related to real data by modifying the optimization in SINDy-like algorithms to allow for more robust model discovery.

SINDy is a dictionary learning approach that represents the nonlinear function $f(X)$ as a linear combination of redundant nonlinear features and applies an iterative sparsity promoting algorithm to find an accurate model and avoid overfitting.  Specifically they solve the overdetermined system $Ax=b$ where $b\approx\dot{X}$, $x$ are the parameters, and $A$ are the features by first solving the least squares problem then removing any terms in the dictionary below a given threshold, and repeating until they acquire a sparse enough solution.  When the data is sufficiently noisy, hard thresholding approaches like this can deactivate important dominant terms early on in the optimization process leaving the algorithm with no mechanism to recover them. As is well known, numerical approximations of the derivative $\dot{u}(t)$ can induce significant noise and this can drastically alter the final extracted model \cite{chartrand2011numerical,chartrand2017numerical}.  One approach is to apply smoothing algorithms to reduce the noise in $b$ \cite{ParamOptimizationDerivatives2020,van2020numerical}, but the use of these methods can introduce biases that mask the true solution \cite{kay2022risk}.

As an alternative approach, we introduce a mechanism to reactivate terms lost in the hard thresholding process with no \emph{a priori} knowledge of the underlying system; we can apply this in many algorithmic settings with hard thresholding. Our motivation follows from simulated annealing algorithms \cite{geman1984stochastic}; we add a ``heat bath" by adding a temperature term to the system and allow the algorithm to ``cool down'' slowly by decreasing said temperature.  In our setting, we reactivate a fraction of the deactivated terms in the thresholding process and slowly let this reactivated fraction reduce to zero.  This allows terms that were eliminated early in the process to be recovered later in the optimization process at the cost of the convergence rate.  In the next section, we will describe our algorithm.  We will compare our algorithm's performance with and without the annealing on various examples highlighting some of the gains that can be made in settings with limited and noisy data, entrainment, and chaos.\\

\begin{figure}
    \includegraphics[width=\linewidth]{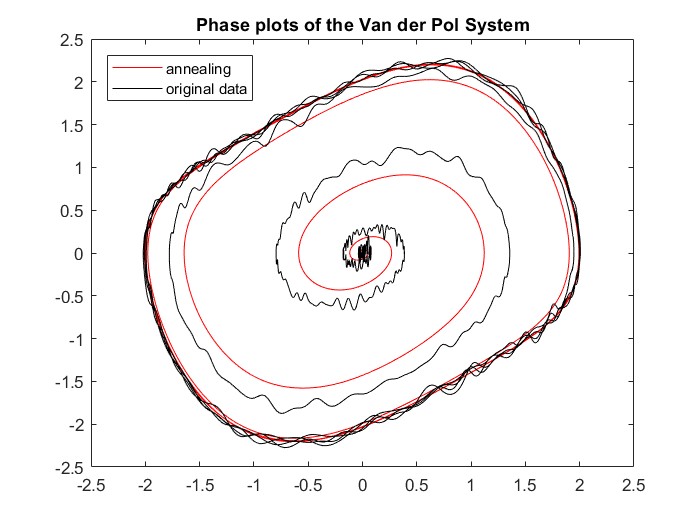}
    \caption{An example plot of the Van der Pol Oscillator with additive zero mean Gaussian noise with variance $0.01$ and a reconstructed solution using SHTreP-A.}
    \label{fig:vdp}
\end{figure}

\noindent\emph{Description of Algorithm} 
Our proposed modifications to hard thresholding schemes can apply in various sparse dictionary learning settings.  First we introduce an alternative to SINDy using a variant of hard thresholding pursuit \cite{foucart2011hard,BLUMENSATH2009265} for the optimization: SINDy with Hard Thresholding Pursuit (SHTreP) and an annealed version (SHTreP-A).  We will apply our annealing approach to \cite{doi:10.1073/pnas.1517384113} yielding SINDy-Anne.  Finally we will compare these four algorithms in the following section.
Given time series data $X$, SHTreP uses a numerical approximation to $\dot{X}$ to find a sparse fit to a library of user-defined functions $\Theta(X)$.  We iteratively seek a sparse coefficient vector $\xi$ by first solving for
\begin{align}\label{gradstep}
    \tilde{\xi}^{k} = H_s(\xi^k + \Theta(X)^T(\dot{X} - \Theta(X)\xi^k))
\end{align}
where $H_s$ is the indicator function on the $s$ largest terms in the vector and $s$ is a hyperparameter that controls the system's sparsity.  This is the hard thresholding process.  We then anneal by reactivating some fraction $0\leq p\leq 1$ of the thresholded terms: this corresponds to considering all entries $\tilde{\xi}_j^{k}=0$ and selecting each of these with probability $p$ for inclusion with the set $S^k=\{1\leq j \leq n: |\hat{\xi}_j^k|>0\}$. Finally we debias by solving the restricted least squares problem on this subset
\[\xi^{k+1}=\mathrm{argmin}_{x\in\mathbb{R}^n:\mathrm{supp}(x)\subseteq S^k} \|\Theta(X)x-b\|_2,\]
giving SHTreP-A.
The initial guess for our algorithm is the least squares solution to $\dot{X} = \Theta(X)\xi$.  The optimization scheme is hard thresholding pursuit from the compressed sensing literature \cite{foucart2011hard}.  If we leave out \eqref{gradstep}, and replace it with a hard thresholding step $\hat{H}_\lambda$ that uses $\hat{S}^k=\{1\leq j \leq n: |\hat{\xi}_j^k|>\lambda\}$ in place of $S_k$, this is SINDy.  In this case, appending random directions like in the above algorithm provides an annealing algorithm SINDy-Anne.  It allows for pruned candidate functions to be reintroduced as the optimization algorithm progresses.

We then iterate this algorithm with an annealing schedule which determines the fraction of terms $p_k$ reactivated at each step $k$.  We will use the annealing schedule given by $\{1, .99,..., .8, .7...,.1, .09, .08,..., .01, 0\}$, meaning we compute one step of sequential hard thresholding then reactivate all of the terms that were removed, simulating a high temperature with lots of movement, then step again and reactivate terms with a probability of 0.99, cooling down, and continue until it no longer reactivates lost terms. The cooling schedule provides another set of hyperparameters for the optimization scheme and can make a significant effect on the obtained minimizer. 
Our cooling schedules will follow the general guidelines in \cite{geman1984stochastic} where the temperature cannot drop too fast, and it must decrease slowly at low temperatures.   After we finish the annealing schedule we let the algorithm run until the termination condition or for around 500 total iterations.\\

\noindent\emph{Comparisons between original and annealed algorithms}
To compare between algorithms with and without annealing, we will typically look at the relative $\ell^2$ and $\ell^1$ errors for ensembles of synthetic examples with progressively larger sampling noise (eg. Fig.~\ref{fig:vdp}):
\begin{itemize}
    \item relative $\ell^2$ error on coefficient vectors $\xi$:
    $$\frac{||\hat{\xi} - \xi||_2}{||\xi||_2} = \sqrt{\frac{\sum_{j=1}^d\sum_{i=1}^n(\hat{\xi}_j^{(i)}-\xi_j^{(i)})^2}{\sum_{j=1}^d\sum_{i=1}^n(\xi_j^{(i)})^2}}$$
    where $\xi_j^{(i)}$ is the coefficient on the $i$th term on the $j$th equation on the original system of equations and $\hat{\xi}_j^{(i)}$ is the $i$th learned coefficient on the $j$th equation in the predicted system.  
    \item relative $\ell^1$ error on coefficients:
    $$\frac{||\hat{\xi} - \xi||_1}{||\xi||_1} = \frac{\sum_{j=1}^d\sum_{i=1}^n|\hat{\xi}_j^{(i)}-\xi_j^{(i)}|}{\sum_{j=1}^d\sum_{i=1}^n|\xi_j^{(i)}|}.$$
\end{itemize}
We will then apply these methods to experimental kinematic data for a tossed foam ball.  The data is limited with a significant amount of structured noise.\\
\begin{figure}
    \includegraphics[width=\linewidth]{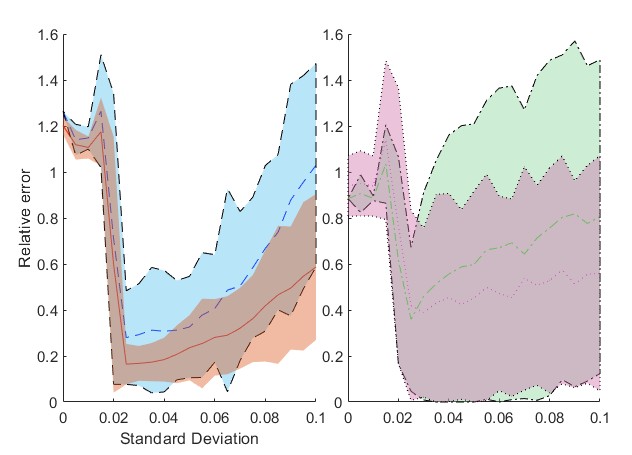}
    \caption{$\ell^1$ coefficient error on 300 different noise profiles for the Lorenz system with the width of the shaded portions being a standard deviation. The blue is SHTreP, the orange is SHTreP-A, the green is standard SINDy, and the purple region is SINDy-Anne. For SINDy, we assume a thresholding parameter for each of the equations: $\lambda_x = 0.4$, $\lambda_y = 0.6$, and $\lambda_z = 0.2$ and for the SHTreP algorithms, we choose a sparsity $s = 15$.}
    \label{fig: lorenz_noise_var}
\end{figure}

\noindent\emph{Lorenz System} First, we consider the popular SciML benchmark Lorenz system \cite{lorenz1963deterministic,gilpin2021chaos} given by
\begin{align*}
    \dot{x} &= 10(y-x)\\
    \dot{y} &= 28x - y - xz\\
    \dot{z} &= xy - 8/3z
\end{align*}
with the initial condition $(x_0,y_0,z_0)^T = (-8,7,27)^T$.  
This is solved forward in time with Runge--Kutta 45 from $t_0 = 0$ to $t_f = 10$ with step size $dt = 0.01$.  We use this solution data and add mean-zero Gaussian noise as our learning data.  At each variance level, we sample 300 instances and record the $\ell^1$ error on the coefficients (Figure \ref{fig: lorenz_noise_var}).  Our library includes 6th degree polynomials and each component equation, $\dot{x},\dot{y},$ and $\dot{z}$ has up to 15 terms. The average error and variance with annealing is lower than without annealing for both algorithms. Finally, note that all the models fixate on incorrect solutions in the absence of sampling noise but actually improve their fit significantly when sampling noise is added. This seems to be common for sparse learning algorithms, but we have not seen it reported in the literature.     

\noindent\emph{FitzHugh--Nagumo}
Next we consider the FitzHugh--Nagumo equations, used in describing excitable systems, given by
\begin{align*}
    \dot{x} &= 0.1 + x + 1/3x^3 - y\\
    \dot{y} &= 0.1(x - y).
\end{align*}
This is also solved with Runge--Kutta 45 with initial condition $(x_0,y_0) = (1,2)$ from $t_0 = 0$ to $t_f = 25$ with time step $0.01$.  For the algorithms, we assume a library of polynomials up to degree 4, with 4 terms in $\dot{x}$ and 3 terms in $\dot{y}$.

\begin{figure}
    \centering
    \includegraphics[width=\linewidth]{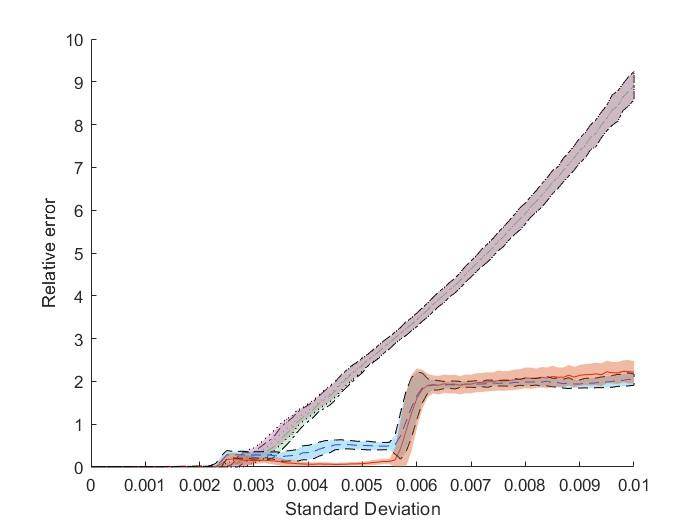}
    \caption{$\ell^1$ error on 300 different noise profiles for the FitzHugh--Nagumo system. The colors are explained in Figure \ref{fig: lorenz_noise_var}. This uses $\lambda = 0.025$ for both $\dot{x}$ and $\dot{y}$ and $s = 6$.}
    \label{fig:fh_ns_l2}
\end{figure}

Figure \ref{fig:fh_ns_l2} reveals similar behavior to the Lorenz system. The annealing again improves the learning in the presence of noise for SHTreP. However, SINDy and SINDy-Anne perform the same, which is expected in systems with low sparsity and a small library.

\begin{figure}
    \centering
    \includegraphics[width=\linewidth]{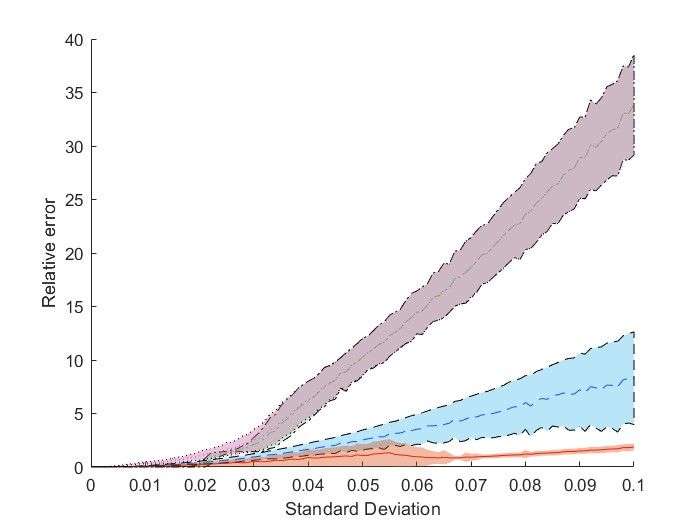}
    \caption{$\ell^1$ error for the Logistic equation starting at $x_0 = 0.01$ to give a sigmoid curve with increasing additive Gaussian noise. This utilizes follows the same color scheme as Figure \ref{fig: lorenz_noise_var}. We use $s = 3$, $\lambda = 0.05$, and a library of degree 4 or less polynomials.}
    \label{fig:log_sig_l1}
\end{figure}

\begin{figure}
    \centering
    \includegraphics[width=\linewidth]{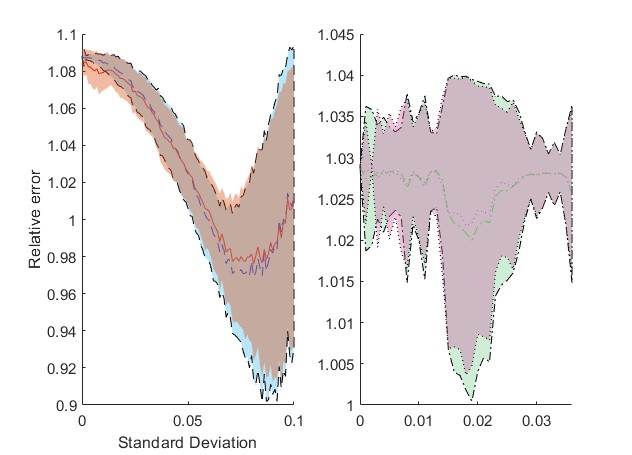}
    \caption{$\ell^1$ error for the Logistic Equation with increasing additive Gaussian noise. This has initial condition $x_0 = 10$ resulting in an exponential decay. The library contains polynomials of degree 10 or less. The hyperparameters are $s = 3$ and $\lambda = 0.001$. The right figure is shortened because the error blows up for the SINDy algorithms.}
    \label{fig:exp_l1}
\end{figure}

\noindent\emph{Logistic}
We now consider the logistic equation:
\begin{align*}
    \dot{x} = rx\left(1-\frac{x}{K}\right)
\end{align*}For our tests, we consider $r = 1, K = 5$, and two different initial conditions:  $x_0 = 0.01$ gives the classic sigmoid curve, and $x_0 = 10$ gives exponential decay to the carrying capacity.  This is solved using the exact solution from $t_0 = 0$, $t_f = 10$, and time step $0.01$. The logistic equation is difficult for sparse learning algorithms becuase it samples so little of phase space.
Figure \ref{fig:log_sig_l1} illustrates that SHTreP-A outperforms all the algorithms as we increase the noise. However, when $x_0 = 10$, we see in Figure \ref{fig:exp_l1} that annealing does not play a significant role for this problem, but all the algorithms generally could not learn with SINDy failing  at a high enough noise.\\

\noindent\emph{Forced Van der Pol Oscillator}
For our final synthetic example, we investigate the Forced Van der Pol Oscillator given by the equation:
\begin{align}
    \ddot{x} + \mu(1-x^2)\dot{x} + x - A\sin(\omega t) = 0.
\end{align}
We take $\mu = 8.53$, $A = 1.2$, and we consider two tests on $\omega$.  This system is nonlinear in the parameters. For the plots of error against noise, we use $\omega = \frac{\pi}{5} + \frac{1}{100}$ which is entrained, and the forcing frequency is overwhelmed by the systems natural frequency.
Then, we conduct a study where $\omega$ is varied from $\frac{\pi}{5}$ to $\pi$. Then, this is solved forward in time starting at initial condition $(x_0,y_0) = (-1,-1)^T$.

In practice, the forcing sinusoids must be included in the library. This is achieved by appending a time-dependent library to the polynomial library where we assume the correct frequency is unknown. We then assume that the frequency can be written as $\omega = \alpha + \varepsilon$ and expand $\sin(\omega t)$ with the angle-sum identity before using a Taylor expansion to get an expression that is quadratic in epsilon. All together, we can approximate the nonlinear forcing with
\begin{align}
    \sin(\omega t) \approx \sin(\alpha t) + \varepsilon t\cos(\alpha t) - \frac{1}{2}\varepsilon^2 t^2\sin(\alpha t). \label{eqn:approx}
\end{align}
This leads to a function library $\Theta(X)$ which contains polynomials of degree 3 or less appended with a new function library 
\begin{align*}
\Phi(\mathbf{t}) = [
\sin(\alpha_1\mathbf{t}) \  ...  \ \sin(\alpha_n\mathbf{t}) \ \mathbf{t}\cos(\alpha_1 \mathbf{t}) \ ... \ \mathbf{t}\cos(\alpha_n\mathbf{t})\\  \mathbf{t}^2\sin(\alpha_1\mathbf{t})  \  ...  \ \mathbf{t}^2\sin(\alpha_n\mathbf{t}) ]
\end{align*}
dependent on the time vector $\mathbf{t} = (t_0,t_1,...,t_f)^T$ and chosen $\alpha_i$ for $i = 1,...,n$. In our tests, we choose $\alpha_i = \frac{\pi}{i}$ for $i = 1,...,6$ near $\alpha=\frac{\pi}{5}$. The augmented function library $[\Theta(X),\Phi(\mathbf{t})]$ contains 44 basis functions (we also include $\cos(\alpha_i\mathbf{t})$ and $\mathbf{t}\sin(\alpha_i\mathbf{t})$ terms), and in practice, the library will contain more as the true $\omega$ remains unknown.

\begin{figure}
    \centering
    \includegraphics[width=\linewidth]{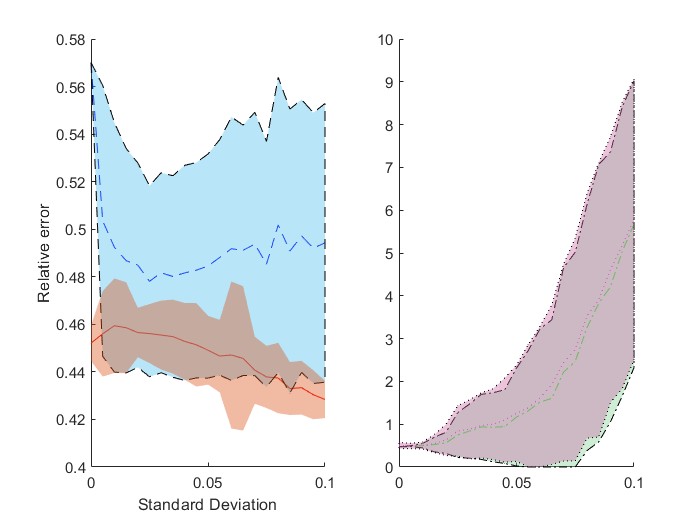}
    \caption{$\ell^1$ error for the Forced Van der Pol Oscillator with increasing additive Gaussian noise, following the conventions from Figure \ref{fig: lorenz_noise_var}. The angular velocity on the forcing term is constant at $\omega = \frac{\pi}{5} + \frac{1}{100}$. The parameters are $s = 10$, $\lambda = 0.005$.}
    \label{fig:fvdp_l1}
\end{figure}

\begin{figure}
    \centering
    \includegraphics[width=\linewidth]{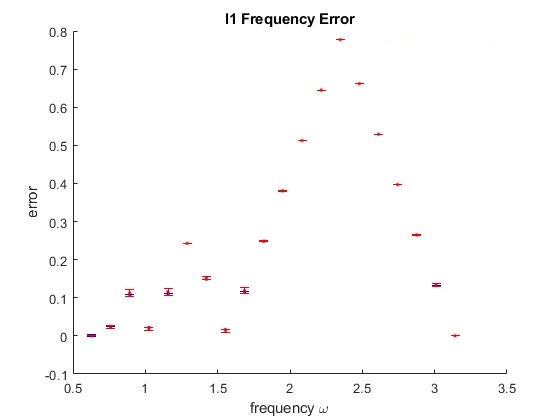}
    \caption{$\ell^1$ error on $\omega$ reconstructed using \ref{eqn:approx} for the Forced Van der Pol Oscillator. The data has additive, zero-mean Gaussian noise with variance $0.05$. The blue represents SHTreP and the red is SHTreP-A}
    \label{fig:freq}
\end{figure}

\begin{figure}
    \centering
    \includegraphics[width=\linewidth]{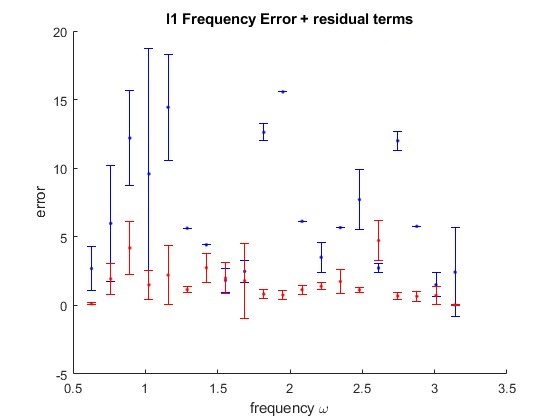}
    \caption{$\ell^1$ error plus the $\ell^1$ norm of the sum of every other time-dependent coefficient on the Forced Van der Pol Oscillator. The data used has additive, zero-mean Gaussian noise with variance $0.05$. The blue is SHTreP and the red is SHTreP-A}
    \label{fig:freq_res}
\end{figure}

Due to the large function library, this is an ideal test of the addition of annealing for use in real-world problems. We conduct the experiment for noise in Figure \ref{fig:fvdp_l1} and notice that SHTreP and SHTreP-A achieve very accurate models despite the large library whereas the SINDy algorithms struggle to learn an accurate model. Then in Figures \ref{fig:freq}, and \ref{fig:freq_res}, we test whether or not the algorithms can recover the correct driving frequency as the sytem varies between chaotic and entrained regimes. The entrained regions can be seen as the tents between frequency that are learned very well.  Furthermore chaotic motion helps the learning process. \\

\begin{figure}
    \centering
    \includegraphics[width=\linewidth]{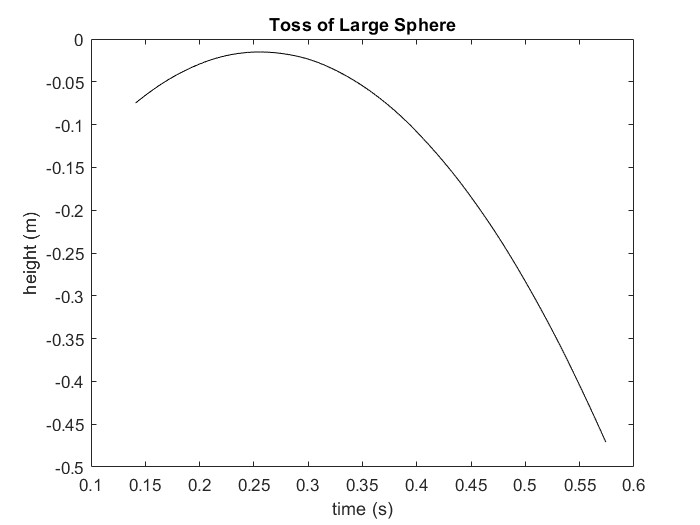}
    \caption{The trajectory of the ball in the y direction.}
    \label{fig:ball}
\end{figure}

\noindent\emph{Projectile Motion}
Model selection algorithms typically struggle with real data; real data often is limited and extremely noisy.
Our final dataset is similar to \cite{de2020discovery} where a foam ball weighing 6.2913 g and of diameter 73.15 mm is tossed up into the air.  The trajectory was collected by applying a particle-tracking algorithm \cite{crocker1996methods} to videos from a high-speed camera (See Figure \ref{fig:ball}).
This data, which is governed by Newton's laws for gravity, was then used to test our algorithms. However, the underlying model drag terms in this model are unknown, so we must use physical intuition to interpret the learned models.  For example, we know $\ddot{y}$ should contain a constant acceleration term of around $9.8$ $\text{m}/\text{s}^2$ and $\ddot{x}$ should not contain any constant force.

We assume that the $x$ and $y$ coordinates decouple and smooth the data with 30 iterations of Matlab's moving average filter, then fit to a library of degree 2 or less polynomials mixed in position and velocity.  In Table~\ref{tab: ball} we see that all the algorithms obtain some of the correct terms, like the constant gravitational acceleration.  However only SHTreP-A finds a reasonable model for the $x$ equation: the SINDy algorithms find a positive $x$ acceleration and SHTreP learns a constant negative driving force.
\begin{table}
    \centering\fontsize{9pt}{8.4pt}
    \begin{tabular}{l l}
        \textbf{Algorithm} & \textbf{Learned Equations}\\
        \midrule
        \multirow{2}{*}{SINDy} & $\ddot{x} = 3.2943x + 5.0572x^2$\\
         & $\ddot{y} = -9.8383 - 60.423y + 15.768y^2 - 3.6329\dot{y}^2$\\
         \hline
        \multirow{2}{*}{SINDy-Anne} & $\ddot{x} = 3.2943x + 5.0572x^2$\\
         & $\ddot{y} = - 9.8383 - 60.423y + 15.768y^2 - 3.6329\dot{y}^2$\\
         \hline
        \multirow{2}{*}{SHTreP} & $\ddot{x} = -4.3883 + 3.4071\dot{x}$\\
          & $\ddot{y} = -8.8419 + 5.3338y + 17.1657y^2$\\
        \hline
        \multirow{2}{*}{SHTreP-A} & $\ddot{x} = -0.6219x - 0.269x\dot{x}$\\
         & $\ddot{y} = -9.062 + 24.0733y^2 - 2.8997y\dot{y}$
    \end{tabular}
    \caption{The learned equations for the projectile motion data over the four algorithms. We assume $s_x = 2$, $s_y = 3$, $\lambda_x = 0.8$, and $\lambda_y = 0.6$.}
    \label{tab: ball}
\end{table}
  \\

\noindent\emph{Conclusion}
This work introduces a novel application of simulated annealing to model identification by changing the hard thresholding process used in sparse learning algorithms.  This allows rejected library terms to be reintroduced leading to improved sparse learning results.  However, there are limitations to this.  On Figure \ref{fig: lorenz_noise_var}, we see that annealing was insufficient to lower error in the low noise case.  Similarly, it was noticed that the annealing sometimes failed when SINDy did not.  We expect that these approaches will be most useful in settings with large function libraries and can be improved by further studying possible annealing schedules  Our preliminary results show promise for these learning algorithms, and can be applied in a variety of settings with hard thresholding.

\begin{acknowledgments}
This work was partially supported by NSF grant $\#2112085$, the Montana State University Undergraduate Scholars Program (USP), and the Hilleman Scholars Programs.  We would also like to thank Benjamin Grodner for providing the projectile motion data.
\end{acknowledgments}

\bibliography{apssamp}
\end{document}